\documentclass[reqno,11pt]{amsart}

\usepackage{latexsym,amssymb,amsthm,graphics,amsmath,amscd}

\usepackage{amssymb}
\usepackage{eucal}
\usepackage{amsmath}
\usepackage{graphicx} 

\theoremstyle{plain}
\newtheorem{theorem}{Theorem}[section]
\newtheorem*{Theorem B}{Theorem B}
\newtheorem*{Theorem A}{Theorem A}

\newtheorem{proposition}{Proposition}[section]
\newtheorem{corollary}{Corollary}[section]
\newtheorem{problem}{Problem}[section]

\newtheorem{example}{Example}[section]
\numberwithin{equation}{section}
\newtheorem {conjecture}{Conjecture}
\theoremstyle{remark}
\newtheorem{remark}{Remark}[section]
 \numberwithin{equation}{section}

\def\<{\left < }
\def\>{\right >}
\def\({\left ( }
\def\){\right )}
\def\e{\eqref}

\def\2{\#_2 }

\begin{document}

\noindent  {\footnotesize Bull. Belg. Math. Soc. Simon Stevin
\vskip-.05in

\noindent 25 (2018), no. 4, pages xxx--xxx}

\vskip.3in 

\title[Two-numbers and their applications - A survey]{Two-numbers and their applications - A survey}

\author[B.-Y. Chen]{Bang-Yen Chen}
\address{Department of Mathematics, Michigan State University, East Lansing, Michigan 48824--1027, U.S.A.}
\email{chenb@msu.edu}

\begin{abstract} The notion of two-numbers of connected Riemannian manifolds was introduced about 35 years ago in [Un invariant geometrique riemannien, C. R. Acad. Sci. Paris Math. 295 (1982), 389--391]  by B.-Y. Chen and T. Nagano. Later, two-numbers have been studied by a number of mathematicians and it was then proved that two-numbers related closely with several important areas in mathematics. 
The main purpose of this article is to survey on two-numbers and their applications. 
At the end of this survey, we present several open problems and conjectures on two-numbers.
\end{abstract}



\dedicatory{In memory of Professor Tadashi Nagano (1930 - 2017)}

\thanks{This is the detailed version of my plenary talk ``Two-numbers and their applications" at the International Conference ``Pure and Applied Differential Geometry - PADGE 2017" held at KU Leuven, Belgium, from August 21 till August 25, 2017.}
\maketitle

\section{two-numbers and maximal antipodal sets}\label{section 1}

The notion of the two-numbers of connected Riemannian manifolds was introduced about 35 years ago  in \cite{CN82}. The primeval concept of  two-numbers is the notion of antipodal points on a circle. 

For a circle $S^1$ in the Euclidean  plane $\mathbb R^2$, the {\it antipodal point} $q$ of a point $p\in S^1$ is the point in $S^1$ which is diametrically opposite to it. 

A geodesic on a Riemannian manifold $M$ is a smooth curve which yields locally the shortest distance between any two nearby points.
Because every closed geodesic in a Riemannian manifold $M$ is isometric to a circle $S^1$,  antipodal points  can be defined for closed geodesics in $M$. More precisely, a point $q$  in a closed geodesic is called an {\it antipodal point} of another point $p$ on the same closed geodesic if the distance $d(p,q)$ between $p$ and $q$ on the two arcs connecting $p$ and $q$ are equal. 

In this article, a closed geodesic in a Riemannian manifold is also called a {\it circle}.
A subset $S$ of a Riemannian manifold $M$ is called an {\it antipodal set} if any two points in $S$ are antipodal on some circle of $M$ connecting them. An antipodal set $A_2M$ in a connected Riemannian manifold $M$ is called a {\it maximal antipodal set} if it doesn't lie in any antipodal set  as a proper subset. The  supremum of the cardinality of all maximal antipodal set of $M$ is called the {\it two-number} of $M$, simply denote by $\#_2M$. 
When a Riemannian manifold $M$ contains no closed geodesics, e.g., a Euclidean space,  we put  
$\#_2M=0$.

Obviously,  we have $\#_2S^n=2$ for any $n$-sphere $S^n$.  On the other hand,  we also have 
\begin{align}\#_2M\geq 2\end{align} 
for every compact connected Riemannian manifold $M$, because compact connected Riemannian manifolds contain at least one close geodesic (cf. \cite{LF51}).

\begin{remark} A maximal antipodal set  $A_2M$ of a Riemannian manifold $M$ is also later known as a {\it  2-set}   or a {\it great set}  for short in some literatures. \end{remark}

\section{Antipodal points in algebraic topology}\label{section 2}

Besides in geometry, the notion of antipodal points plays some significant roles in many areas of mathematics, physics, and applied sciences. 

The following result from algebraic topology is well-known. 

\vskip0.1in
\noindent {\bf The Borsuk--Ulam Antipodal Theorem.} {\it Every continuous function from an $n$-sphere $S^n$ into the Euclidean $k$-space $\mathbb E^k$ with $k\leq n$ maps some pair of antipodal points to the same point. 
In other word, if  $f: S^n\to \mathbb E^k$  is continuous, then there exists an $x\in S^n$ such that
$f(-x)=f(x)$.}

\vskip0.1in
Clearly, Borsuk-Ulam's theorem fails for $k>n$, because $S^n$ embedded in $\mathbb E^{n+1}$. The Borsuk-Ulam theorem has numerous applications; range from combinatorics to differential equations and even economics.

\vskip.05in
For $n=1$, Borsuk--Ulam's theorem implies that  that {\it at any moment there always exist a pair of opposite points on the earth's equator with the same temperature.}  This assumes the temperature varies continuously.

For $n=2$, Borsuk--Ulam's  theorem implies that {\it  there is always a pair of antipodal points on the Earth's surface with equal temperatures and equal barometric pressures},  at any moment.

The Borsuk--Ulam theorem was conjectured by S. Ulam at the Scottish Caf\'e in Lvov, Ukraine. 
\noindent Ulam's conjecture was solved in 1933 by K. Borsuk \cite{B33}.  It turned out that
the result had been proved three years before  in \cite{LS30} by L. Lusternik and L. Schnirelmann (cf. \cite[footnote, page 190]{B33}).
 Since then many alternative proofs and also many extensions of Borsuk--Ulam's theorem have been found as collected by H. Steinlein in his survey article \cite{S85} including a list of 457 publications involving various generalizations of Borsuk-Ulam's theorem.

Another important theorem from algebraic topology with the same flavor is the following.

\vskip0.1in
\noindent {\bf Brouwer's Fixed-Point Theorem.} {\it Every continuous function  $f: B\to B$ 
    from the unit $n$-ball $B$ into itself has a fixed point.} 

\vskip0.1in

Brouwer's fixed-point theorem has numerous applications to many fields as well. 
For instance,  Brouwer's fixed-point theorem and its extension by S. Kakutani in \cite{IKa41} (extending  Brouwer's fixed-point theorem to set-valued functions) played a central role in the proof of existence of general equilibrium in Market Economies as developed in the 1950s by economics Nobel prize winners K. Arrow (1972) and G. Debreu (1983).

\begin{remark} In fact,  Borsuk-Ulam's antipodal theorem implies Brouwer's fixed point theorem, see \cite{Su97}.
\end{remark}

\section{Symmetric spaces}\label{section 3}

The class of Riemannian manifolds with parallel Riemannian curvature tensor, i.e., $\nabla R=0$, was first introduced independently by P. A. Shirokov  in 1925 and by H. Levy in 1926. This class is known today as the class of {locally symmetric Riemannian spaces} (see, e.g., \cite[page 292]{C00}). 

It was \'E. Cartan who noticed in 1926 that irreducible spaces of this type are separated into ten large classes each of which depends on one or two arbitrary
integers, and in addition there exist twelve special classes corresponding to the
exceptional simple groups. Based on this, \'E. Cartan created his theory of symmetric Riemannian spaces in his famous  papers ``Sur une classe remarquable d'espaces de Riemann'' in 1926 \cite{C26}.

Symmetric spaces are the most beautiful and important Riemannian manifolds. Such spaces arise in a wide variety of situations in both mathematics and physics.  
 This class of spaces contains many prominent examples which are of great importance for various branches of mathematics, like compact Lie groups, Grassmannians and bounded symmetric domains.  Symmetric spaces are also important objects of study in representation theory, harmonic analysis as well as in differential geometry.

An isometry $s$ of a Riemannian manifold is called an {\it involutive} if $s^{2}=id$. A Riemannian manifold $M$ is called a {\it symmetric space} if for each point $x\in M$ there is an involutive isometry $s_{x}$ such that $x$ is an isolated fixed point of $s_{x}$; the involutive isometry $s_{x}\ne id$ is called the {\it symmetry at $x$}. A Hermitian symmetric space is a Hermitian manifold $M$ such that every point of $M$ admits a symmetry preserving the Hermitian structure of $M$.

Let $M$ be a symmetric space. Denote by $G$ the closure of the group of isometries on $M$ generated by $\{s_p:p\in M\}$ in the compact-open topology. Then $G$ is a Lie group which acts transitively on the symmetric space; hence the typical {\it isotropy subgroup} $K$, say at a point $o\in M$, is compact and $M=G/K$.  
Thus symmetric spaces are  homogeneous spaces as well. 
From the point of view of Lie theory, a symmetric space is the quotient $G/K$ of Lie group G by a Lie subgroup $K$, where the Lie algebra  ${\mathfrak {k}}$ of $K$ is also required to be the 
(+1)-eigenspace of an involution of the Lie algebra ${\mathfrak {g}}$ of $G$. 

In this paper, we will use standard symbols as in Helgason's book \cite{H} to denote symmetric spaces mostly. Here are a few minor exceptions.
More specifically than $AI, AI(n)$ denote $SU(n)/ SO(n), AI\hskip-.02in I(n) :=SU(2n)/ Sp(n)$,
etc. Let $G_d({\mathbb  R}^n)$, $G_d({\mathbb  C}^n)$ and $G_d({\mathbb  H}^n)$ denote the Grassmannians of $d$-dimen\-sional subspaces in the real, complex and quaternion vector spaces (or modules), respectively.
The standard notations for the exceptional spaces such as $G_2, F_4,E_6,. . . ,$ $GI,. . . ,EI\hskip-.02in X$ denote the simply-connected spaces where we write $GI$ for $G_2/ SO(4)$.
We will denote by $M^*$ the  bottom space, i.e., the adjoint space in \cite{H}, of the space $M$. 
 
For a symmetric space $M$, the dimension of a maximal flat totally geodesic submanifold of $M$ is a well-defined natural number; called the {\it rank} of $M$ and denoted by $rk(M)$. 
Clearly, the rank of a symmetric space is at least one. It is well-known that the class of rank one compact symmetric spaces consists of $n$-sphere $S^n$, a projective space ${\mathbb F}P^n ({\mathbb F}={\mathbb  R}, {\mathbb C}, {\mathbb H})$ and  the 16-dimensional Cayley plane $FI\hskip-.02in I={\mathcal O}P^2$ with ${\mathcal O}$ being the  Cayley  algebra. 

Obviously, every complete totally geodesic submanifold of a symmetric space is also a symmetric space. 
 It follows from the equation of Gauss that 
\begin{align} rk(B)\leq rk(M)\end{align} 
for every complete totally geodesic submanifold $B$ of a symmetric space $M$ (cf. \cite{C73}).

\section{$(M_+,M_-)$--theory}\label{section 4} 

In this section, we provide a  brief introduction of the $(M_+,M_-)$--theory for compact symmetric spaces introduced by B.-Y.  Chen and  T. Nagano in \cite{CN77,CN78,CN88,Na88,Na92}. Our approach to compact symmetric spaces based on antipodal points and fixed point sets of compact symmetric spaces. Hence our approach to compact symmetric spaces is quite different from other approaches done by \'E. Cartan and others.

In fact, our approach to compact symmetric spaces plays the Key Roles in the determination of two-numbers of compact symmetric spaces; including  the determination of 2-ranks of compact Lie groups (see \cite{C87,C13}).

 Let $o$ be a point of a compact symmetric space $M=G/K$. A connected component of the fixed point set $F(s_{o},M)\smallsetminus \{o\}$ of the symmetry $s_{o}$ is called a {\it polar} of $o$. We denote it by $M_{+}$ or $M_{+}(p)$ if $M_{+}$ contains a point $p$.  
 
 We have the following useful propositions from \cite{CN78} (also \cite[page 15]{C87}).

\begin{proposition}\label{P:4.1} Let $M=G/K$ be a compact symmetric space. Then for each antipodal point $p$ of $o\in M$, the isotropy subgroup $K$ at $o$ acts transitively on the polar $M_+(p)$. Moreover, we have $K(p)=M_+(p)$ and $K(p)$ is connected. Consequently, we have $M_+(p)=K/K_p$, where $K_p$ is given by $\{k\in K: k(p)=p\}$.
\end{proposition}
 
When a polar consists of a single point, we call it a {\it pole}. 
 
 \begin{proposition}\label{P.4.2} Under the hypothesis of Proposition \ref{P:4.1}, the normal space to $M_+(p)$ at $p$ in $M$ is the tangent space of a connected complete totally geodesic submanifold  $M_-(p)$. Thus we have
  \begin{align}\dim M_{+}(p)+\dim M_{-}(p)=\dim M.\end{align}
\end{proposition}

 \begin{proposition}\label{P:4.3} For each antipodal point $p$ of $o$ in a compact symmetric space $M$, we have
 {\rm (1)}  $rk (M_-(p))=rk (M)$ and
{\rm (2)} $M_-(p)$ is a connected component of the fixed point set $F(s_{p}\circ s_{o},M)$ of $s_{p}\circ s_{o}$ through $p$.
\end{proposition}

Polars and meridians of a compact symmetric space are totally geodesic submanifolds; they are thus compact symmetric spaces as well. Both polars and meridians have been determined for every compact connected irreducible Riemannian symmetric space (see  \cite{C87,CN88,Na88,Na92}). One of the most important properties of polars and meridians is that $M$ is determined (globally) by any pair of $(M_{+}(p),M_{-}(p))$. 

 Let $o$ be a point of a compact connected Riemannian symmetric space $M$. If there exists a pole $p$ of $o$ in $M$, then the set consisting of the midpoints of all geodesics  joining $o$ and $p$ is called the {\it centrosome} of $\{o,p\}$; denoted by $C(o,p)$. 
  Every connected component of  $C(o,p)$ is again a totally geodesic submanifold of $M$. 
   Centrosomes play important roles in topology as well. For instance, centrosomes have been used by J. M. Burns in \cite{Bu92} to compute homotopy in many compact symmetric spaces (see subsection 6.4).  
 
 The following result from \cite[page 277]{CN88} characterizes poles in compact symmetric spaces (see also \cite{C87.2}).

\begin{proposition} \label{P:4.4}
 The following six conditions are equivalent to each other for two distinct points $o,p$ of a connected compact symmetric space $M=G_{M}/K_{G}$.

\begin{itemize} 

\item[(i)] $p$ is a pole of $o\in M$; 
  
\item[(ii)] $s_{p}=s_{o}$; 
  
\item[(iii)] $\{p\}$  is a polar of $o\in M$; 
  
\item[(iv)] there is a double covering totally geodesic immersion $\pi=\pi_{\{o,p\}} :M\to M''$ with $\pi(p)=\pi(o)$;

\item[(v)] $p$  is a point in the orbit $F(\sigma,G_{M})(o)$  of the group $F(\sigma,G_{M})$ through $o$, where $\sigma = {\rm ad}(s_{o})$;

\item[(vi)] the isotropy subgroup of $SG_{M}$ at $p$ is that, $SK_{G}$ $($of $SG_{M}$ at $o)$, where $SG_{M}$ is the group generated by $G_{M}$ and the symmetries; $SG_{M}/G_{M}$ is a group of order $\leq 2$.
\end{itemize}  
\end{proposition}
  
 For a compact symmetric space $M$, the {\it Cartan quadratic morphism} $$Q=Q_{o}: M\to G_{M}$$ carries a point $x\in M$ into $s_{x} s_{o}\in G_{M}$. The Cartan  quadratic morphism is a $G_{M}$-equivariant morphism which is an immersion.
  
 We have the following result  for centrosomes from \cite[pages 279-280]{CN88}.
  
\begin{proposition} \label{P:4.5}   The following five conditions are equivalent to each other for two distinct points $o,q$ of a connected compact symmetric space $M$.

\begin{itemize} 

\item[(i)] $s_{o}s_{q}=s_{q}s_{o}$; 
  
\item[(ii)] $Q(q)^{2}=1_{G_{M}}$, where $Q=Q_{o}$ is Cartan  quadratic morphism; 
  
\item[(iii)] either $s_{o}$ fixes $q$ or $q$ is a point in the centrosome $C(o,p)$ for some pole $p$ of $o$; 
  
\item[(iv)] either $s_{o}(q)=q$ or $s_{o}(q)=\gamma (q)$ for the covering transformation $\gamma$ for some pole $p=\gamma(o)$ of $o$;

\item[(v)] either $s_{o}(q)=q$ or there is a double covering morphism $\pi:M\to M''$ such that $s_{o''}$ fixes $q''$, where $o''=\pi(o)$ and $q''=\pi(q)$.
\end{itemize}  
\end{proposition}

 The polars, meridians and centrosomes play very important roles in the study of compact connected symmetric spaces as well as of compact connected Lie groups.  In particular, 
Propositions \ref{P:4.4} and \ref{P:4.5} play important roles for the study of maximal antipodal sets and two-numbers of compact symmetric spaces as well as the 2-ranks of compact Lie groups.
  
Professor T. Nagano and M. Sumi (= M. S. Tanaka) proved in \cite{NM88} that, for a compact symmetric space $M$, the root system $R(M_-)$ of a meridian $M_-\ne M$ is obtained from the Dynkin diagram of the root system $R(M)$ of $M$ and they also found all maximal totally geodesic spheres in $SU(n)$ by means of the $(M_+,M_-)$-method.

 \section{Two-numbers of compact symmetric spaces}\label{section 5}

\subsection{Antipodal set in terms of symmetries}\label{section 5.1}

For compact symmetric spaces, we have the following result form  \cite[page 275]{CN88}.

\begin{proposition} \label{P:5.1} For any compact symmetric space $M$,  the two-number $\#_2M$  is equal to the maximal possible cardinality $\#(A_2M)$ of a subset $A_2M$ of $M$ such that the point symmetry $s_x$ fixes every point of $A_2M$ for every $x\in A_2M$. 
\end{proposition}

For compact symmetric spaces, Proposition \ref{P:5.1} can be regarded as an alternative definition of two-numbers.

The following inequality from  \cite[page 276]{CN88} is an important result for two-numbers.

\begin{theorem} \label{T:5.1} $\2 M-1$ does not  exceed the sum of the two-numbers of all the polars of a point in a compact symmetric space $M$, that is
\begin{align}\label{5.1} \2 M\leq 1+\sum \2 M_+.\end{align}
\end{theorem}

\begin{remark} The equality \e{5.1} holds in many cases (such as the groups
$Sp(n)$ and $O(n)$ and the hermitian symmetric spaces) and does not in the other
cases (such as the adjoint group $SU(8)^*$ of $SU(8)$).
\end{remark}

\subsection{Two-numbers as an obstruction to totally geodesic embeddings}\label{section 5.2}

If $B$ is a totally geodesic submanifold of a Riemannian manifold $M$, then every geodesic of $B$ is also a geodesic of $M$. Thus we have \begin{align}\#_2B\leq \#_2M.\end{align} Consequently,
the invariant $\#_2M$ is a geometric obstruction to the existence of a totally geodesic embeddings $f: M\to N$, since the existence of $f$
clearly implies the inequality \begin{align}\label{5.3}\#_2M\leq \#_2N.\end{align} 

\vskip.05in
For instance,  although the real projective $n$-space $RP^n$ can be topologically embedded in $S^m$ for sufficient high $m$, but  inequality \e{5.3} implies that every real projective $k$-space $RP^k$ with $k\geq 2$ cannot totally geodesic embedded into $S^m$ regardless of codimension. This is simply due to the fact:
$$\#_2 RP^k=k+1 > 2=\#_2 S^m.$$

Similarly, while the complex Grassmann manifold $G_2 ({\mathbb C}^4)$ of the 2-dimen\-sional subspaces of the complex vector space ${\mathbb C}^4$ is obviously embedded into $G_3({\mathbb C}^6)$ as a totally submanifold, the space $G_2 ({\mathbb C}^4)^*$ which one obtains by identifying every member of $G_2 ({\mathbb C}^4)$ with its orthogonal complement in ${\mathbb C}^4$, however, cannot be totally geodetically embedded into $G_3({\mathbb C}^6)^*$
simply due to
$$\#_2 G_2 ({\mathbb C}^4)^*=15>12 =G_3 ({\mathbb C}^6)^*$$ according to \e{5.1}.

\subsection{Two-numbers of dot products}\label{section 5.3}

Suppose a finite group $\Gamma$ is acting on two spaces $M$ and $N$ freely as automorphism groups. Then $\Gamma$ acts on the product space $M\times N$ freely. And the orbit space $(M\times N)/\Gamma$ is called the {\it dot product of M and
N} (with respect to $\Gamma$) and denoted by $M\cdot N$. In most cases of our study, $\Gamma$ will be the group of order two acting on $M$ and $N$ as the covering transformation groups for double
covering morphisms in the sequel. $\Gamma$ will not be mentioned in that case, if $\Gamma$ is
obvious or if $\Gamma$ need not be specified. 

\begin{example} $SO(4) = S^3\cdot S^3$, $U(n) = T\cdot SU(n)$ and $GI$ has the only polar
$S^2\cdot S^2$. Here $\Gamma$ for $U(n)$ is the center of $SU(n)$, a cyclic group of order $n$.
\end{example}

The following result from \cite[page 281]{CN88} provides simple links between dot products, centrosomes and two-numbers of compact symmetric spaces.

\begin{theorem}\label{T:5.2} The dot product for double coverings $M\to M''$ and $N\to N''$ has the following properties: 

 \begin{itemize}
 \item[{\rm (a)}]  $\#_2M\leq \#_2(M\cdot N)$;
 \item[{\rm (b)}] $\#_2(M\cdot N)\leq \#_2(F(s_o,M)\cdot F(s_p,N))+\#_2(CM\cdot CN)$, where $CM$ and $CN$  are the centrosomes for the point $o$ of $M$ and its pole and for $p$of $N$ and its pole;
  \item[{\rm (c)}] $\frac{1}{2}(\#_2 M)(\#_2  N)\leq \#_2(M\cdot N)\leq 2 (\#_2M'')(\#_2 N'')$.
 \end{itemize}
\end{theorem}

The following corollaries are easy consequences of Theorem \ref{T:5.2}.

\begin{corollary} For a compact symmetric space $M$, we have 
$$ \#_2 (S^n\cdot M) \leq  \#_2 M + \#_2 (S^{n-1}\cdot CM),$$ 
where $CM$ is a centrosome of $M$.
\end{corollary}

\begin{corollary} For a compact symmetric space $M$, we have 
$$ \#_2 M\leq \2 (S^1\cdot M) \leq  \#_2 M + \#_2  CM\leq 2\2 M.$$
\end{corollary}

\begin{corollary} We have 
$ \#_2 (S^m\cdot S^n) =2(n+1)$ if $m\geq n$.
\end{corollary}

\subsection{Two-numbers of irreducible compact symmetric spaces}\label{section 5.4}
In \cite{CN88},  Nagano and I have determined  the two-numbers and the maximal antipodal sets for most irreducible compact symmetric spaces (see also Appendices II and III of \cite[pages 67--71]{C87}). However, the antipodal sets of {\it oriented} real Grassmannians have not been discussed in  \cite{CN88}. 

On the other hand,  H. Tasaki had described the maximal antipodal sets of oriented real Grassmannians in \cite{Tas13,Tas15}.

\section{Two-numbers and Topology}\label{section 6}

The two-numbers link closely with topology. In this section, we shall present a number of results in this respect.

\subsection{Two-numbers and Euler numbers}\label{section 6.1}

First, we present some links between two-numbers and Euler numbers of compact symmetric spaces from \cite[Theorem 4.1]{CN88}.

\begin{theorem} \label{T:6.1}  For  any compact symmetric space $M$, we have
\begin{align}\label{6.1}\#_2M\geq\chi(M),  \end{align}
where $\chi(M)$ denotes the Euler number of $M$.
\end{theorem}

The proof of this theorem based on a result of Hopf-Samelson on rank of $G$ and $K$ of a homogeneous space $G/K$ with $\chi(G/K)$,  a result of H. Hopf on fixed point sets as well as the $(M_+,M_-)$-theory of symmetric spaces.

For any hermitian compact symmetric space of semisimple type, Professor Nagano and I proved in \cite[Theorem 4.3]{CN88} the following.

\begin{theorem}\label{T:6.2}   For  any compact hermitian symmetric space $M$ of semi-simple type, we have
\begin{align}\label{6.2} \#_2M=\chi(M)=1+\sum \2 M_+.\end{align}
\end{theorem}
The proof of this theorem based heavily on the $(M_+,M_{-})$--theory as well as the Lefschetz
fixed point theorem in the version of Atiyah and Singer \cite{AS68}. 

\vskip.1in
An immediate consequence of this theorem is the following.

\begin{corollary} For every complete totally geodesic hermitian subspace $B$ of a semi-simple hermitian symmetric space $M$, we have
\begin{align}\label{6.3}\chi(M)\geq\chi(B).  \end{align}
\end{corollary}

\subsection{Two-numbers and homology}\label{section 6.2}

A symmetric R-space is a special kind of compact symmetric space for which several characterizations were known. 
Originally, Professor T. Nagano  introduced the notion of a symmetric R-space  in \cite{Na65}   as a compact symmetric space $M$ which admits a Lie transformation group $P$ which is non-compact and contains $I_0(M)$ as a subgroup. For example, a sphere is a symmetric R-space. The
$I(M)$ of $M = S^n$ coincides with the natural action of $O(n + 1)$ on $S^n$. Thus the
action of $SO(n+1)$ on $S^n$ coincides with $I_0(M)$.
 
 In the theory of algebraic groups, a {\it Borel subgroup} of an algebraic group $G$ is a maximal Zariski closed and connected solvable algebraic subgroup.  Subgroups between a Borel subgroup $B$ and the ambient group $G$ are called {\it parabolic subgroups}.  Working over algebraically closed fields, the Borel subgroups turn out to be the minimal parabolic subgroups in this sense. Thus $B$ is a Borel subgroup when the homogeneous space $G/B$ is a complete variety which is as large as possible.

M. Takeuchi used the terminology {\it symmetric R-space} for the first time in \cite{Ta65}.  He gave a cell decomposition of an R-space in \cite{Ta65}, which is a kind of generalization of a symmetric R-space. Here, by a {\it R-space} we mean $M = G/U$ where $G$ is a connected real semisimple Lie group without center and $U$ is a parabolic subgroup of $G$.
   
   A compact symmetric space $M$ is said to have a {\it cubic lattice} if a maximal torus of $M$ is isometric to the quotient of ${\mathbb E}^r$ by a lattice of ${\mathbb E}^r$ generated by an orthogonal basis of the same length. 
     O. Loos \cite{Lo85} gave another intrinsic characterization of symmetric R-spaces among all compact symmetric spaces with the property that the unit lattice of the maximal torus of the compact symmetric space (with respect to a canonical metric) is a cubic lattice. Loos' proof was based on the correspondence between the symmetric R-spaces and compact Jordan triple systems. 

In  \cite{Fe80}, D. Ferus characterized symmetric R-spaces as compact symmetric submanifolds of Euclidean spaces. He also proved that compact extrinsically symmetric submanifolds are orbits of isotropy representations of symmetric space of compact type (or non-compact type).

 \vskip.1in
 The class of symmetric $R$-spaces includes:
\begin{itemize}
  \item[{\rm (a)}] all hermitian symmetric spaces of compact type

 \item[{\rm (b)}]  Grassmann manifolds $O(p+q)/O(p)\times O(q), Sp(p+q)/Sp(p)\times Sp(q)$

 \item[{\rm (c)}]  the classical groups $SO(m),\,U(m),\,Sp(m)$,

 \item[{\rm (d)}]  $U(2m)/Sp(m),\, U(m)/O(m)$,

 \item[{\rm (e)}]  $(SO(p+1)\times  SO(q+1))/S(O(p)\times O(q))$, where $S(O(p)\times O(q))$ is  the subgroup  of $SO(p+1)\times SO(q+1)$ consisting of matrices of the form:
$$\begin{pmatrix} \epsilon& 0 & & \\ 0 & A&&\\ &&\epsilon& 0\\
&&0&B \end{pmatrix},\quad \epsilon=\pm1,\quad A\in
O(p),\quad B\in O(q),$$
 (This $R$-space is cover twice by $S^p\times S^q$),

 \item[{\rm (f)}]  the Cayley projective plane $FII={\mathcal O}P^2$, and 

 \item[{\rm (g)}]  the three exceptional spaces $E_6/Spin(10)\times T, E_7/E_6\times T,$ and $E_6/F_4.$
\end{itemize}

For any symmetric $R$-space,  M. Takeuchi proved in \cite{Ta89} the following. 

\begin{theorem}\label{T:6.4}  For  any symmetric $R$-space $M$, we have
\begin{align}\label{6.6} \#_{2}M=\sum_{i\geq 0} b_i(M,\mathbb Z_2),  \end{align}
where $b_i (M,\mathbb Z_{2})$ is the i-th Betti number of $M$ with coefficients in $\mathbb Z_2$.
\end{theorem}

This theorem was proved by applying a result of Chen-Nagano from \cite{CN88} in conjunction with an earlier result of M. Takeuchi from \cite{Ta65}.

\subsection{Bott's periodicity theorem for homotopy of classical groups}\label{section 6.3}

The most famous work of R. Bott is his  periodicity theorem  which describes periodicity in the homotopy groups of classical groups (cf. \cite{Bo70}).

 Bott's original results may be succinctly summarized as 
 
 \begin{theorem}\label{T:6.5} The homotopy groups of the classical groups are periodic:
$$ \pi_i(U)=\pi_{i+2}(U),\; \pi_i(O)=\pi_{i+4}(Sp),\; \pi_i(Sp)=\pi_{i+4}(O)$$
for $i=0,1,\cdots$,
where $U$ is the direct limit defined by $U=\cup_{k=1}^\infty U(k)$ and similarly for $O$ and $Sp$.
\end{theorem}
 
 The second and third of these isomorphisms given in Theorem \ref{T:6.5} yield  the {\it  8-fold periodicity}: 
$$\pi_i(O)=\pi_{i+8}(O),\; \pi_i(Sp)=\pi_{i+8}(Sp)$$
for $i=0,1,\cdots.$

Bott's original proof of his periodicity theorem  is differential geometric
in its nature (see \cite{Bo59}). His proof relies on the observation that in a compact Riemannian
symmetric space $M$ one can choose two points $p$ and $q$ in ``special position''
such that the connected components of the space of shortest geodesics in $M$
joining $p$ and $q$ are again compact symmetric spaces. 

Put $M=M_0$ and let $M_1$ be one of the resulting connected components. This construction can be repeated inductively. 

Let $M$ be a compact symmetric space. The {\it index of a geodesic} $\gamma$ in $M$ from $p$ to $q$ is the number of conjugate points of $p$ counted with their multiplicities, in the open geodesic segment from $p$ to $q$. 

We denote the space of shortest geodesics from $p$ to $q$ by the {\it symbol $\Omega^d$}, which relates closely with the notion of centrosomes.

The proof of Bott's periodicity theorem in \cite{Bo59} relied on the following result.

\vskip.1in
\noindent {\bf Rott's theorem.} {\it If $\Omega^d$ is a topological manifold and if every non-shortest geodesic from $p$ to $q$ has index greater than or equal to $\lambda_0$, then the $(i+1)$-th homotopy groups of $M$ satisfies
\begin{align} \pi_{i+1}(M)\cong \pi_i(\Omega^d)\end{align}
for $i<\lambda_0-1$.}

\subsection{Applications to homotopy}\label{section 6.4}

The $(M_+,M_-)$--theory as well as centrosomes play significant roles  in computing homotopy of compact symmetric spaces.   

In fact, J. M. Burns mentioned in the introduction of his paper \cite{Bu92} that  how our $(M_+, M_-)$--theory (in particular, centrosomes)  in  conjunction with Bott's theorem given above can be used   to compute the homotopy of compact symmetric spaces.
In \cite{Bu92}, J. M. Burns carried out   the computation of the homotopy in compact symmetric spaces of types:  $AI, AI\hskip-.02in I$ and $ CI$.   
 By applying the same method, he also computed in \cite{Bu92} the homotopy in the exceptional symmetric  spaces:  $ EI\hskip-.02in I\hskip-.02in I-EI\hskip-.02in X$, $F_4,\; FI$ and $ FI\hskip-.02in I.$

The homotopy groups $\pi_i(EV\hskip-.02in  I\hskip-.02in I)$ with $i\leq 16$ were computed in 2012 by P. Quast \cite{{Qu12}}. The first eight homotopy groups of $EV\hskip-.02in  I\hskip-.02in I$  have been already determined by J. M. Burns in \cite{Bu92}.

\begin{remark} It is known that  the homotopy group $\pi_i(M)$ of a connected manifold $M$ is the same as the homotopy group $\pi_i(\tilde M)$ of its universal cover $\tilde M$  for $i\geq 2$.
Also, the homotopy groups $\pi_9, \pi_{10}$ and $\pi_{14}$ of $EV\hskip-.02in  I\hskip-.02in I$ can also be directly read off from the long exact sequence of
homotopy groups of coset spaces together with the corresponding homotopy groups of $E_6$ and $E_7$ that can be found in \cite[page 363]{MT91}.

\end{remark}

\subsection{Two-numbers and covering maps}\label{section 6.5}

 In \cite{CN88}, Professor T. Nagano and I  also discovered some simple links between two-numbers and covering maps between compact symmetric spaces.
   
For double coverings we have \cite[page 280]{CN88}:
  
\begin{theorem} \label{T:6.6} If $M$ is a double covering of $M''$,  then $\#_{2}M\leq 2\cdot \#_2(M'').$
 \end{theorem}
 
\begin{remark} Th inequality in Theorem \ref{T:6.6} is sharp, since the equality case holds for $M=SO(2m)$ with $m>2$.
\end{remark}

 For $k$-fold coverings with odd $k$, \cite[page 278]{CN88}:
 
\begin{theorem} Let $\phi :M\to N$ is a $k$-fold covering  between compact symmetric spaces.  If $k$ is odd then $\#_{2}M=\#_2(N).$
 \end{theorem}

\subsection{Antipodal sets, 2-numbers and Borsuk-Ulam's type theorems}\label{section 6.7}

A continuous function $f: M\to {\mathbb R}$ of a compact symmetric space $M=G/K$ into the real line is called {\it isotropic} if $f$ is invariant under the action of the isotropic subgroup $K$. 

Recently, I obtained the following Borsuk-Ulam's type results involving antipodal sets and two-numbers of compact symmetric spaces in \cite{C17}.

\begin{theorem}\label{T:6.8} Let  $f: M\to {\mathbb R}$ be a continuous isotropic function from a compact symmetric space into the real line. Then  $f$ carries  some maximal antipodal set of $M$ to the same point  in ${\mathbb R}$, whenever $M$ is one of the following spaces: Spheres; the projective spaces ${\mathbb F}P^n\, ({\mathbb F}={\mathbb  R}, {\mathbb C}, {\mathbb H})$;  the Cayley plane $FI\hskip-.02in I$; the exceptional spaces $EIV; EIV^*;GI$; and the exceptional Lie group $G_2$.
Consequently, the function $f: M\to {\mathbb R}$ carries a maximal  antipodal set with $\#_2M$ elements to the same point  in ${\mathbb R}$.
\end{theorem}

All of the compact symmetric spaces $M$ mentioned in the list of Theorem \ref{T:6.8} have single polar for a fixed point  $o\in M$. Next, we discuss the case in which the point $o\in M$ has multiple polars. 
Suppose $o\in M$ has multiple polars, say $M_+^1, M_+^2,\ldots,M_+^k$. Then we denote by $\hat M_+$ a polar of $o$ with the maximal two-number among all polars of $o$. 

We also have the following results from \cite{C17}.

\begin{theorem}\label{T:6.9}  Let  $f: M\to {\mathbb R}$ be a continuous isotropic function of a compact symmetric space $M$ into the real line ${\mathbb R}$. If $M$ admits more than one polar,
then  $f$ carries an antipodal set $S_M$  of $M$ consisting of  $1+\#_2\hat M_+$ points of $M$ to the same point  in ${\mathbb R}$.
\end{theorem}

The following corollaries are easy consequences of Theorem \ref{T:6.8} (see \cite{C17} for more details).

\begin{corollary}\label{C:6.2}  If  $f: E_8\to {\mathbb R}$ is a continuous isotropic function of the exceptional Lie group $E_8$ into the real line ${\mathbb R}$,
then  $f$ carries an antipodal set  of $E_8$ with  $392$ elements to the same point  in ${\mathbb R}$.
\end{corollary}

\begin{corollary}\label{C:6.3}  If  $f: FI\to {\mathbb R}$ is a continuous isotropic function of the compact symmetric space $FI$ into the real line ${\mathbb R}$,
then  $f$ carries an antipodal set  of $FI$ with  $24$ elements to the same point  in ${\mathbb R}$.
\end{corollary}

We provide a simple example to illustrate that the isotropic condition on $f$ in Theorem \ref{T:6.8} is essential.

\begin{example} {\rm Consider the real projective plane ${\mathbb R}P^2$. Then there is a double covering $\pi:S^2 \to {\mathbb R}P^2$. Without loss of generality, we may assume that $S^2$ is the unit sphere in $\mathbb E^3$ centered at the origin of $\mathbb E^3$. 
For each real-valued function $f:{\mathbb R}P^2\to {\mathbb R}$, the lift $\hat f:S^2\to {\mathbb R}$ of $f$ is an even function via the double covering $\pi$, i.e., $\hat f(-{\bf x})=\hat f({\bf x}),\;\;  \forall {\bf x}=(x,y,z)\in S^2.$

Conversely, for every given real-valued even function $g:S^2 \to \mathbb R$ of $S^2$, $g$ induces a real-valued function $\check{g}:{\mathbb R}P^2\to \mathbb R$ of ${\mathbb R}P^2$.

If we choose $g:S^2\to \mathbb R$ to be  $g({\bf x})=(x-y)^2$,
then $g$ induces a real-valued function $\check{g}:{\mathbb R}P^2\to \mathbb R$ which does not carry any maximal antipodal set of ${\mathbb R}P^2$ to the same point in $\mathbb R$.}
\end{example}

\section{2-rank of compact connected Lie groups}

\subsection{2-ranks of Borel and Serre}\label{section 7.1}

A. Borel and J.-P. Serre considered in \cite{BS53} the 2-rank, denoted by $r_2G$, of a compact connected Lie group $G$.  By definition the 2-rank of a compact connect Lie group $G$ is the maximal possible rank of the elementary 2-subgroup of $G$. 

In \cite{BS53} Borel and Serre proved the following 2 results.
\vskip.05in

\begin{itemize}
  \item[{\rm (1)}] 
 $rk(G)\leq r_2(G)\leq 2rk(G)$ with $rk(G)=rank (G)$ and 
 
 \item[{\rm (2)}] $G$ has (topological) 2-torsion if $rk(G)<r_2(G)$. 
\end{itemize}
\vskip.05in

 Borel and Serre were able to determine  the 2-rank of  simply-connected simple Lie groups $SO(n), Sp(n), U(n), G_{2}$ and $F_{4}$  in \cite{BS53}. Furthermore, they proved that the exceptional Lie groups $G_{2}$, $F_{4}$ and $E_{8}$ have 2-torsion. 
On the other hand, they also mentioned in  \cite[page 139]{BS53} that they were unable to determine the 2-rank for the exceptional simple Lie groups $E_{6}$ and $E_{7}$.

\subsection{2-rank and commutative algebra}\label{section 7.2}

After Borel and Serre's paper \cite{BS53},  2-ranks of compact Lie groups $G$ have been investigated by a number of mathematicians. In particular, it had been shown that the 2-ranks of compact Lie groups have several important links with commutative algebra.
Here, we  provide merely two such links.

\vskip.1in
(a) Let $F$ be either a field or the rational integer ring $\mathbb Z$. Let 
$A=\sum_{i\geq 0} A_{i}$
 be a graded commutative $F$-algebra in sense of Milnor-Moore \cite{MM65}. If $A$ is connected, then it admits a unique augmentation $\varepsilon :A\to F$.   
 
 Put  $\bar A={\rm Ker}\, \varepsilon$. The $\bar A$ is called the {\it augmentation ideal} of $A$.
 A sequence of elements $\{x_{1},\ldots,x_{n}\in \bar A\}$ in $\bar A$ is called a {\it simple system of generators} if $\{x_{1}^{\epsilon_{1}}\cdots x_{n}^{\epsilon_{n}}: \epsilon_{i}=0 \hbox{ or }1\}$ is a module base of $A$.

 Let $G$ be a compact connected Lie group. Denote by $s(G)$ the number of generators of a simple system of the $\mathbb Z_{2}$-cohomology $H^{*}(G,{\mathbb Z}_{2})$ of  $G$. 

A. Kono discovered in \cite{Ko77} the following relationship between $s(G)$,  $r_{2}G$ and the $\mathbb Z_{2}$-cohomology $H^{*}(G,{\mathbb Z}_{2})$.

 \begin{theorem}\label{T:7.1} Let $G$ be a connected compact Lie group. Then
  the following three conditions are equivalent: 
  
  \begin{itemize}
  \item[{\rm (1)}] $s(G)\leq r_{2}G$; 
  
  \item[{\rm (2)}] $s(G)= r_{2}G$; 
  
  \item[{\rm (3)}] $H^{*}(G,{\mathbb Z}_{2})$ is generated by universally transgressive elements. 
\end{itemize}
\end{theorem}

To prove Theorem \ref{T:7.1}, A. Kono had applied May's spectral sequence \cite{Ma65},  Eilenberg-Moore's spectral sequence \cite{EM62} as well as Quillen's result in \cite{Qu71.1}. 
In \cite{Ko77}, Kano also described some properties of compact Lie groups satisfying condition (3) in Theorem \ref{T:7.1} and provided applications.

\vskip.1in

(b) In commutative algebra, {Krull's dimension} of a ring $R$ is the supremum of the number of strict inclusions in a chain of prime ideals. More precisely, one says that a strict chain of inclusions of prime ideals of the form: 
 $${\mathfrak p}_{0}\subsetneq {\mathfrak p}_{1}\subsetneq \cdots \subsetneq {\mathfrak p}_{n}$$
 is of length $n$. That is, it is counting the number of strict inclusions. 
   
Given a prime ideal ${\mathfrak p}\subset R$, we define the {\it height} of ${\mathfrak p}$ to be  the supremum of the set 
$\{n\in {\mathbb N}: {\mathfrak p} \hbox{ is the supremum of a strict chain of length $n$}\}.$
Then the Krull dimension is the supremum of the heights of all of its primes.
 
 Let $G$ be a compact Lie group and put $H_{G}^{*}=H^{*}(BG; {\mathbb Z}_{2}),$ where $BG$ is a classifying space for $G$. Let $N_{G}^{*}\subset H_{G}^{*}$ denote the ideal of nilpotent elements.  Then it is known that $H_{G}^{*}/N_{G}^{*}=H_{G}^{\#}$  is a finitely generated commutative algebra.  
 
 In \cite{Qu71.1}, D. Quillen studied the relationship between the finitely generated commutative algebra $H_{G}^{\#}$  and the  structure of the Lie group $G$. In particular, he proved that, under some suitable assumptions, the Krull dimension of $H_{G}^{\#}$ is equal to the 2-rank of $G$. 
 
 D. Quillen proved his result by calculating the mod 2 cohomology ring of extra special 2-groups. 
  Quillen's result provided an affirmative answer to a conjecture of M. F.  Atiyah \cite{At61} and also of a conjecture of R. G. Swan  \cite{Sw71}.

 \section{Applications of two-numbers to group theory}
 
 If $G$ is a connected compact Lie group, then by assigning $s_{x}(y)=xy^{-1}x$ to every point $x\in G$, we have $s^{2}_{x}=id_{G}$ to each point $x$. Therefore $G$ becomes a compact symmetric space with respect to a bi-invariant Riemannian metric. 
 
\subsection{Links between two-numbers and 2-ranks}\label{section 8.1}
  The following simple link between the 2-rank and the two-number of a connected compact Lie group $G$ was established in \cite[Proposition 5]{CN82}.
 
 \begin{theorem}\label{T:8.1} Let $G$ be a connected compact Lie group. Then we have
 \begin{equation} \#_{2}G=2^{r_{2}G}.\end{equation}
 \end{theorem}
 
 For products of two compact Lie groups, we have the following result from \cite[Lemma 1.7]{CN88}.
 
\begin{theorem}\label{T:8.2} Let $G_{1}$ and $G_{2}$ be  connected compact Lie groups. Then 
 \begin{equation}\label{8.2} \#_{2}(G_{1}\times G_{2})=2^{r_{2}G_{1}+r_{2}G_{2}}.\end{equation}
\end{theorem}

Based on these two theorems and  $(M_+,M_-)$--theory, we were able to determine the 2-ranks of  compact connected Lie groups in \cite[pages 289-293]{CN88}.  
Consequently, we settled the problem of Borel-Serre for determination of 2-ranks of  connected compact simple Lie groups as follows.

 \subsection{2-ranks of classical groups}\label{section 8.2} 
 
 For classical groups we have:

\begin{theorem}\label{T:8.3} Let $U(n)/{\mathbb Z}\mu$  by the quotient group of the unitary group $U(n)$ by the cyclic normal subgroup ${\mathbb Z}\mu$ of order $\mu$. Then we have 
  \begin{equation} r_{2}(U(n)/{\mathbb Z}\mu)= 
 \begin{cases} n+1 & \text{if $\mu$ is even and $n=2$ or $4$;}\\ n & \text{otherwise.}\end{cases} \end{equation}
 \end{theorem}

 \begin{theorem}\label{T:8.4} For $SU(n)/{\mathbb Z}\mu$, we have 
 \begin{equation} r_{2}(SU(n)/{\mathbb Z}\mu)= 
 \begin{cases} n+1 & \text{for $(n,\mu)=(4,2)$;}
 \\ n  & \text{for $(n,\mu)=(2,2)$ or $(4,4)$;}
 \\ n-1 & \text{for the other cases.} \end{cases} \end{equation}
 \end{theorem}

 \begin{theorem}\label{T:8.5} One has $r_{2}(SO(n))=n-1$ and, for $SO(n)^{*}$, we have
 \begin{equation} r_{2}(SO(n)^{*})= 
 \begin{cases} 4 & \text{for $n=4$;}\\ n-2 & \text{for $n$ even $>4$}.\end{cases} \end{equation}
 \end{theorem}

 \begin{theorem}\label{T:8.6} Let $O(n)^{*}=O(n)/\{\pm 1\}$. We have
 
 {\rm (a)} $r_{2}(O(n))= n$;
 
 {\rm (b)} $r_{2}(O(n)^{*})$ is $n$ if $n$ is 2 or 4, while it is $n-1$ otherwise.
 \end{theorem}

 \begin{theorem}\label{T:8.7} One has $r_{2}(Sp(n))=n$, and, for $Sp(n)^{*}$, we have
  \begin{equation} r_{2}(Sp(n)^{*})= \begin{cases} n+2 & \text{for $n=2$ or $4$}\\ n+1 & \text{otherwise}.
\end{cases} \end{equation}
Thus we also have
 \begin{equation}r_{2}(Sp(n)^{*})= r_{2}(U(n)/{\mathbb Z}_{2}) + 1 \end{equation} for every $n$. 
\end{theorem}

 \subsection{2-ranks of spinors,  semi-spinors and $Pin(n)$}\label{section 8.3}
 
Next, we consider the spinor $Spin(n)$ and its related groups. Recall that  $Spin(n)$ is a subset of the Clifford algebra $Cl(n)$ which is generated over $\mathbb R$ by
the vectors ${\bf e}_{i}$ in the fixed orthonormal basis of $\mathbb R^{n}$; subject to the conditions
${\bf e}_{i}{\bf e}_{j}=-{\bf e}_{j}{\bf e}_{i}$ and $ {\bf e}_{i}{\bf e}_{i}= -1,\, i\ne j.$

 \vskip.1in

For $Spin(n)$ we have the following two results.

 \begin{theorem}\label{T:8.8} We have
 \begin{equation}\notag
  r_{2}(Spin(n)) =\begin{cases} r+1 & \text{if  $\,n \equiv -1,0$ or {\rm 1 (mod 8)}}\\ 
 r &\text{otherwise},\end{cases}
 \end{equation} where $r$ is the rank of $Spin(n)$, $r = [\frac{n}{2}]$. 
 \end{theorem}

 \begin{theorem}\label{T:8.9} {\rm (PERIODICITY)} For $n\geq 0$, One has
$$ r_{2}(Spin(n+8))= r_{2}(Spin(n))+4$$
 
\end{theorem}

The group $Pin(n)$ was introduced in \cite{ABS64} by M. F. Atiyah, R. Bott and A. Shapiro while they studied Clifford modules. 
$Pin(n)$ is a group in the Clifford algebra $Cl(n)$ and it double covers $O(n)$ and whose connected component $Spin(n)$ double covers $SO(n)$.

For $Pin(n)$ we have:

 \begin{theorem}\label{T:8.10} For $Pin(n)$ we have $r_{2}(Pin(n))= r_{2}(Spin(n + 1))$, $n\geq 0$.
 \end{theorem}

For the semi-spinor group $SO(4m)^{\#}=Spin(4m)/\{1, e_{((4m))}\}$, we have:

 \begin{theorem}\label{T:8.11}   We have
 \begin{equation}\notag
  r_{2}(SO(4m)^{\#}) =\begin{cases} 3 & \text{if  $m=1$}\\  6 &\text{if $m=2$,}\\ r+1 &\text{if $m$ is even $>2$, }\\
 r &\text{if $m$ is odd $>1$},\end{cases}
 \end{equation}
 where $r$ is the rank $2m$ of $SO(4m)^{\#}$.  
 \end{theorem}

\begin{remark}The 2-rank of $Spin(16)$ and of $SO(16)^{\#}$ had been obtained independently by J. F. Adams in \cite{Ad87}.  However, the method of his proof is completely different from ours given in \cite{CN88}.
\end{remark}

 \subsection{2-ranks of exceptional groups}\label{section 8.4}

We also able to determine the 2-ranks of Exceptional groups as follows.

\begin{theorem}\label{T:8.12} One has $$r_{2}G_{2}=3,\;\;  r_{2} F_{4}=5,\;\;  r_{2}E_{6}=6,\;\; r_{2}E_{7}=7
,\;\; r_{2}E_{8}=9$$ for the simply-connected exceptional simple Lie groups.\end{theorem}

\begin{remark} $r_2 G_2=3$ and $r_2 F_4=5$ are due to A. Borel and J. P. Serre.
\end{remark}

We also able to determine the rank of the {bottom space} $E^{*}_{6}$.

 \begin{theorem}\label{T:10.13} One has $r_{2}E^{*}_{6}=6$. \end{theorem}

\section{Two-numbers and algebraic geometry}

\subsection{Two-numbers and projective ranks}\label{section 9.1}

The {\it projective rank, $Pr(M)$,} of a compact hermitian symmetric space 
     $M$ is  the maximal complex dimension of  totally geodesic complex projective spaces $N$ of $M$ defined by A. Fauntleroy in  \cite{Fa93}.  
     
     A subset $\mathcal E$ of a maximal antipodal set of a compact symmetric space $M$ is called {\it equidistant} if there is a real number $a>0$ such that $d_M(x,y)=a$ for any two points $x,y\in \mathcal E$. 
     
     Put
\begin{align} \gamma=\gamma_M=\min \{d_M(x,y): x,y\in {\mathcal E}, x\ne y\}\end{align}
and let $A_\gamma\subset \mathcal E$ be an equidistance set (for the distance $\gamma$) of maximal cardinality. Let $\mu=\# (A_\gamma)$ denote the cardinality of $A_\gamma$.

 C. U. S\'anchez and A. Guinta  proved in \cite{SG02} the following.
 
 \begin{theorem} \label{T:9.1} Let $M$ be a compact connected irreducible hermitian symmetric space.
 
 \begin{enumerate}  
\item[{\rm (a)}] If $M\ne CI(n)$, then $Pr(M)=\mu-1$.
  
\item[{\rm (b)}] If $M=CI(n)$ and $A_{\sqrt{2}\gamma}$ is an equidistant set of maximal cardinality (for the distance $\sqrt{2}\gamma(M))$ and $\mu$ is its cardinality, then $Pr(M)=\mu-1$.
  \end{enumerate}
  \end{theorem}
  
 Theorem \ref{T:9.1}  provides the following simple link between the two-number and the projective rank for compact irreducible hermitian symmetric spaces. 
 \begin{theorem} \label{T:9.2} $Pr(M)\cdot rk(M)\leq \#_2(M)$ for every compact irreducible hermitian symmetric space $M$.
  \end{theorem}

\subsection{Two-numbers and arithmetic distances}\label{section 9.2}

W.-L. Chow introduced in \cite{Ch49} the notion of {\it arithmetic distance ``$d$''} for every classical hermitian symmetric spaces $M$.
In case that the classical hermitian symmetric space $M$ is the complex Grassmannian $G_p({\mathbb C}^n)$,  W.-L. Chow's arithmetic distance $d(V,W)$ between $V,W\in G_p({\mathbb C}^n)$ is defined by
\begin{align}d(V,W)=\dim_{\mathbb C} V/V \cap W,\end{align}
i.e., the codimension of their intersection $V \cap W$ in $V$ (or in $W$). Equivalently, the arithmetic distance $d(V,W)$ is the smallest integer $t$ such that there exists a finite set $\{U_i\}_{1\leq i \leq t+1}$ of linear subspaces in $G_p({\mathbb C}^n)$ such that 
\begin{enumerate}  
\item[{\rm (a)}] $U_1=V$ and $U_{t+1}=W$;
  
\item[{\rm (b)}]  $d(U_i,U_{i+1})=1$ for $1\leq i \leq t$. 
    \end{enumerate}

Chow proved  in \cite{Ch49} that the {\it $d$-preserving transformations of $M$ are either holomorphic or anti-holomorphic, provided $rk (M)>1$}
 
 \vskip.1in
A {\it Helgason sphere} in a compact symmetric space  is a maximal-dimensional totally geodesic sphere with maximal sectional curvature. The dimension of a Helgason sphere is thus equal to one more than the multiplicity of the highest restricted root.
  
Using Helgason spheres, S. Peterson \cite{Pe87} extended  Chow's arithmetic distance  to an arithmetic distance $d$ to {arbitrary irreducible compact symmetric spaces}.  In this general case, $d(x,y)$ satisfies
 $d(x,y)\leq j$  
 if the two points $x,y\in M$ are joined by a chain of $j$  Helgason spheres.

 More precisely,  put $d(x,y)=0$ if $x=y$. And put $d(x,y)=1$ if $x\ne y$ and $x,y$ lie in a Helgason sphere.
 Otherwise, $d(x,y)$ is defined to be the smallest  $j$ such that there is a chain of $j$ Helgason spheres joining $x$ and $y$  (instead of a chain of linear subspaces in $G_p({\mathbb C}^n)$). 
 
  If $M$ is hermitian and irreducible, the holomorphic transformations permute the Helgason spheres and hence they preserve the new arithmetic distance obviously. Also, it equals Chow's arithmetic distance if $M=G_p({\mathbb C}^n)$.

  Let $M=G/K$ be a compact Riemannian symmetric space where $G$ denotes the identity component of the isometry group of $M$.   T. Nagano studied and classified in \cite{Na65} the {\it geometric transformation groups} of compact symmetric spaces. 
  Roughly speaking, they are ``larger groups'' $L$ that act on $M$ such that 
  
\begin{enumerate}  
\item[{\rm (i)}] $G\subset L$; 
  
\item[{\rm (ii)}]  $L$ is a Lie transformation group acting effectively on $M$; 
  
\item[{\rm (iii)}]  $L$ preserves the symmetric structure of $M$; and 
  
\item[{\rm (iv)}]  $L$ is simple. 
  \end{enumerate}

By using our approach to compact symmetric spaces and a formulation of Radon's duality, S. Peterson proved the following rigidity theorem  in \cite{Pe87}.

\begin{theorem}\label{T:9.3} Let $M$ be a $G_d({\mathbb F}^n),\, {\mathbb F}={\mathbb R}, {\mathbb C}, {\mathbb H}$, or $M=AI(n)$ with $\dim M \geq 3$ and let
$L=\{\varphi: M\to M,$  $\varphi$ a diffeomorphism  preserving
the arithmetic distance $d\}.$
Then $L=L'$,  where $L'$ be the geometric transformation group of $M$.
\end{theorem}

\section{Further applications of two-numbers}\label{section 10}

\subsection{Real forms of hermitian symmetric spaces}\label{section 10.1}  

Let $M$ be a hermitian symmetric space of compact type and let $\tau$ be an  involutive anti-holomorphic isometry of $M$. The fixed point set $$F(\tau,M) = \{x \in M : \tau(x) = x\}$$ is called a {\it real form} of $M$,    which is connected and a totally geodesic Lagrangian submanifold. 
    For instance, let $$M = CP^1 = S^2,\;\;  L = RP^1 = S^1.$$ Then $M$ is a hermitian symmetric space of compact type and $L$ is a real form of $M$. 
     Any two distinct great circles in $S^2$ intersect at two points which are antipodal to each other. 
     The intersection is an antipodal set in $L$ as well as in $M$.
 
M. Takeuchi proved the following result in \cite{Ta84}.
  
\begin{theorem} \label{T:10.1} Every real form of a hermitian symmetric space of compact type is a symmetric $R$-space. 
Conversely, every symmetric $R$-space
is realized as a real form of a hermitian symmetric space of compact type. 

Moreover, the correspondence is one-to-one.\end{theorem} 

In \cite{TT12}, M. S. Tanaka and H. Tasaki studied the intersection of two real forms in a hermitian symmetric space of compact type. In particular, they proved that the intersection of two real forms in a Hermitian symmetric space of compact type is an antipodal set if the intersection is discrete.
More precisely, they established  the following link between real forms of hermitian symmetric spaces  and two-numbers in \cite{TT12}.   

\begin{theorem}  \label{T:10.2} Let $L_1,L_2$ be real forms of a hermitian symmetric space 
of compact type whose intersection is discrete. Then $L_1\cap L_2$ is an
antipodal set in $L_1$ and $L_2$.
Furthermore,   if $L_1$ and $L_2$ are congruent, then $L_1\cap L_2$
is a 2-set. Thus $\#(L_1\cap  L_2) = \#_2L_1 = \#_2L_2.$
\end{theorem} 

Further results in this respect see  \cite{Ik15,TT15,TT15.1}.

\subsection{Fixed point set of a holomorphic isometry}\label{section 10.2}

O. Ikawa, M. S. Tanaka and H. Tasaki showed in \cite{ITT15}  a necessary and sufficient condition
that the fixed point set of a holomorphic isometry of a Hermitian symmetric space
of compact type is discrete and proved that the discrete fixed point set is an antipodal
set. They also proved a necessary and sufficient condition that the intersection of two
real forms in a Hermitian symmetric space of compact type is discrete and consider
some relations between the intersection of two real forms and the fixed point set of
a certain holomorphic isometry by the use of the symmetric triads introduced by
O. Ikawa in \cite{Ik11}.

\subsection{Applications to Lagrangian Floer homology}\label{section 10.3} 

Let $(M,\omega)$ be a symplectic manifold, i.e., $M$ is a manifold with a closed nondegenerate 2-form $\omega$. Let $L$ be a Lagrangian submanifold in $M$, i.e.,
$\dim_{\mathbb R} L = \frac{1}{2} \dim_{\mathbb R} M$ and $\omega$ vanishes on $L$.

The symplectic {\it Floer homology} is a homology theory associated to a symplectic manifold and a nondegenerate symplectomorphism of it. If the symplectomorphism is Hamiltonian, the homology arises from studying the symplectic action functional on the (universal cover of the) free loop space of a symplectic manifold. 

Symplectic Floer homology is invariant under Hamiltonian isotropy of the symplectomorphism.
Let $Hamilt(M,\omega)$ denote the set of all Hamiltonian diffeomorphisms of $M$.

A. Floer defined in \cite{Fl88}  the homology when $\pi_2(M,L_i)=0$, $i=0,1$, and
proved that it is isomorphic to the singular homology group $H_*(L_0,\mathbb Z_2)$ of $L_0$ in the
case where $L_0$ is Hamiltonian isotopic to $L_1$. 
As a result,  Floer solved affirmatively the so called {\it Arnold conjecture} for Lagrangian intersections in that case. 

In \cite{Gi89}, A. Givental posed the following conjecture which generalized the results of Floer and himself.

\vskip.1in
\noindent {\bf Conjecture (Arnold-Givental).} {\it Let $(M,\omega)$ be a symplectic manifold and
$\tau: M\to M$ be an anti-symplectic involution of $M$. Assume that the fixed point set
$L=F(M,\tau)$ is not empty and compact. Then for any $\phi\in Hamilt(M,\omega)$ such that the
Lagrangian submanifold $L$ and its image $\phi(L)$ intersect transversally, the inequality
\begin{align}\label{10.1} \#(L\cap \phi(L))\geq \ b(L,\mathbb Z_2)\end{align}
holds, where $b(L,\mathbb Z_2)=\sum_{i\geq 0} b_i(L,\mathbb Z_2)$.}
\vskip.1in

Recall that M. Takeuchi proved that if $L$ is a symmetric $R$-space, then one has $\#_2 L=b(L,\mathbb Z_2)$ (cf. Theorem \ref{T:6.4}).
H. Iriyeh, T. Sakai and H. Tasaki  computed in  \cite{IST13}
 the Lagrangian Floer homology $HF(L_0,L_1;{\mathbb Z}_2)$  of a pair of real forms 
$(L_0,L_1)$ in a monotone hermitian symmetric space $M$ of compact type in the case where $L_0$ is not necessarily congruent to $L_1$. In particular, they obtained a generalization of the Arnold-Givental inequality \e{10.1} in the case where $M$ is irreducible. 

As an application, H. Iriyeh, T. Sakai and H. Tasaki proved the following.

\begin{theorem} Every totally geodesic Lagrangian sphere in the complex hyperquadric is globally volume minimizing under Hamiltonian deformations.
\end{theorem}

\subsection{Two-numbers and tight Lagrangian submanifolds}\label{section 10.4} 

The  complex hyperquadric $Q_{n}(\mathbb C)$ is holomorphically isometric to the compact hermitian symmetric space $\tilde G_{2}({\mathbb R}^{n+2})$. 
 
    For every natural number $k$ with $0\leq k\leq n$, the complex hyperquadric $Q_{n}(\mathbb C)$ admits certain  real forms $S^{k,n-k}$ defined by $S^{k,n-k}=(S^{k}\times S^{n-k})/{\mathbb Z}_{2}.$
 
H. Tasaki  proved in \cite{Ta10} the following.

\begin{theorem} \label{T:10.4}   Let $k,\ell$ be integers with $0\leq k\leq \ell \leq [\frac{n}{2}]$, and let $L_{1}$ and $L_{2}$ be real forms  which are congruent to $S^{k,n-k}$ and $S^{\ell,n-\ell}$, respectively.  If $L_{1}$ and $L_{2}$ intersect transversally, then $L_{1}\cap L_{2}$ is a 2-set of $L_{1}$ and an antipodal set of $L_{2}$. Moreover, if $k=\ell=[\frac{n}{2}]$,  then   the two-numbers satisfy $$\#_2(\tilde G_{2}({\mathbb R}^{n+2}))=\#(L_{1}\cap L_{2}).$$
\end{theorem}  

 By applying Theorem \ref{T:10.4}, H. Tasaki obtained the following.

\begin{theorem} \label{T:10.5}  Any real form of the oriented real Grassmannian $\tilde G_{2}({\mathbb R}^{n+2})$ is a globally tight Lagrangian submanifold.\end{theorem}

\subsection{Application to convexity}\label{section 10.5} 

Let $M$ be a Riemannian manifold and let $\tau$ be an involutive isometry  of $M$. A connected component of the fixed point set of $\tau$ with positive dimension is called a {\it reflective submanifold}, which is a totally geodesic submanifold of $M$.
   A connected submanifold $S$ of $M$ is called  (geodetically) {\it convex} if any shortest geodesic segment in $S$ is still shortest in $M$.
 
The $(M_+,M_-)$--theory  was also used in \cite{QT12}  by P. Quast and M. S. Tanaka to prove the following result.

\begin{theorem}\label{T:10.6} Every reflective submanifold of a symmetric $R$-space is convex.\end{theorem}

\subsection{Applications to stability}\label{section 10.6} 

By applying the $(M_+,M_-)$--theory for compact symmetric spaces,  the author, P.-F. Leung and T. Nagano  introduced in 1980  an algorithm in \cite{CLN80} to determine the stability of totally geodesic submanifolds of compact symmetric spaces. 

In 1987, Y. Ohnita reformulated and improved  the above algorithm in \cite{Oh87} to include the formulas for the index, the nullity and the Killing nullity of a compact totally geodesic submanifold
in a compact symmetric space.
 In the same paper, Ohnita proved that every Helgason sphere in any compact symmetric space is stable. Subsequently, many results about the stability of totally geodesic submanifolds in compact symmetric spaces were obtained by applying our algorithm and Ohnita's reformulation.

In 1990, K. Mashimo and H. Tasaki proved in \cite{MT90} that every closed Lie subgroup of Dykin index 1 in every compact simple connected Lie group is stable. 

K. Mashimo also
determined in \cite{Ma92} all of unstable Cartan embeddings of compact symmetric spaces. 
Moreover, by applying the algorithm there are results on the stability of symmetric $R$-spaces in Hermitian symmetric spaces obtained by M. Takeuchi \cite{Ta84} and also the stability of polars and meridians by M. S. Tanaka in \cite{Ta95}.

Furthermore, T. Kimura and M. S. Tanaka determined in \cite{KT09} the stability of the maximal totally geodesic submanifolds in compact symmetric spaces of rank two by applying the algorithm.

\subsection{Two-numbers and algebraic coding theory}\label{section 10.7}
Coding theory is the study of the properties of codes and their respective fitness for specific applications. 
The main purpose of codes is to be able to recover the original content of a transmitted message by correcting errors that have entered the message during transmission. This capability is useful in maintaining the integrity of communication systems, computer networks, compact disk recording, etc.

A $p$-group $H$ is called {\it extra special} if its center ${\mathbb Z}$ is cyclic of order $p$, and the quotient $H/{\mathbb Z}$ is a non-trivial elementary abelian $p$-group.

J. A. Wood explored in \cite{Wo89} the equivalence between the diagonal extra-special 2-group 
 of the connected compact Lie group $Spin(n)$ and the self-orthogonal linear binary codes of algebraic coding theory. Our results in subsection \ref{section 8.3} (Theorems 8.8 and 8.9) were mentioned and been used in his work in \cite{Wo89}.

\subsection{An application to theory of designs}\label{section 10.8}
The theory of designs is the part of combinatorial mathematics that deals with the existence, construction and properties of systems of finite sets whose arrangements satisfy generalized concepts of balance and/or symmetry.

In 1973, P. Delsarte \cite{De73} unified the theories of codes and designs on association schemes  and gave the upper bounds for codes and the lower bound for designs by applying linear programming for polynomials associated with metric or cometric association schemes. 

After Delsarte's work, the theory of spherical designs was introduced by Delsarte-Goethals-Seidel in \cite{DGS77} as an analogy of Delsarte technique. The essential tool in their works is the addition formula for polynomials; polynomials associated with metric or cometric association schemes, or the Gegenbauer polynomials with spheres. 

In general, the theory of designs can be given on the Delsarte spaces  or the polynomial spaces, which are metric spaces with ``good'' polynomials, such as the polynomials associated with metric or cometric association schemes or  Gegenbauer's polynomials. 

Compact symmetric spaces of rank one are natural and significant examples of the Delsarte spaces or the polynomial spaces for continuous metric spaces. The theory of designs on compact symmetric spaces of rank one was studied in details by S. G. Hoggar in \cite{Ho82}.

For compact symmetric spaces of higher rank, H. Kurihara and T. Okuda proved  in \cite{KO13}  a characterization of maximal antipodal sets of complex Grassmannians in term of certain designs (more precisely, ${\mathcal E}\cup {\mathcal F}$-designs) with the smallest cardinalities.

\section{Extensions of two-numbers and applications}\label{section 11}

\subsection{Holomorphic two-numbers of compact hermitian symmetric spaces}\label{section 11.1}

C. U. S\'anchez  defined  {\it holomorphic two-number} $\#_2^H(M)$ for a compact connected hermitian manifold $M$ in \cite{Sa97} as the maximal possible cardinality of a subset $A_2$ such that for every pair of points $x$ and $y$ of $A_2$,  there exists a totally geodesic complex 
curve of genus 0 in $M$ on which $x$ and $y$ are antipodal to each other.

C. U. S\'anchez proved the following.

\begin{theorem} \label{T:11.1} $\#_{2}^H M=\#_2(M)\,$ for every compact hermitian symmetric space. 
\end{theorem}

By combining Theorem \ref{T:11.1} with our equality $\#_2M=\chi(M)$ from Theorem \ref{T:6.2}, one has the following.
 
\begin{corollary}  \cite{Sa97} For every compact hermitian symmetric space $M$, one has $\#_2^H(M)$ $=\chi(M)$.
\end{corollary}

  \subsection{Index numbers for real flag manifolds}\label{section 11.2}
 
A {\it real flag manifold} (also called a {\it R-space}) is a homogeneous space of the form $G/P$, where $G$ is a real semisimple Lie group without compact factors and $P$ is a parabolic
subgroup. Here $G/P\equiv K/K \cap P$ which is a $K$-orbit on $P$.
 
Every complex flag manifold may be considered an $R$-space. In fact, let $U$ be a compact connected semisimple centerless Lie group and $\frak u$ its Lie algebra. The
{\it complex flag manifolds} of $U$ are the orbits of the adjoint action of $U$ on $\frak u$. 

Take
$M= Ad(U)Y$ for $Y \ne 0$ in $U$. Let ${\mathfrak g}= {\mathfrak u}_c= {\mathfrak u}+i {\mathfrak u}$. Then there is a Cartan decomposition of the realization ${\mathfrak g}_{\mathbb R}$  of ${\mathfrak g}$ and one may consider $M$ as the orbit of $iY$  in $i{\mathfrak u}$ by the adjoint action of $U$.
 
 In \cite{Sa97}, C. U. S\'anchez proved the following.

\begin{proposition}\label{P:11.1} Let $M$ be a real flag manifold. Then there exists a complex flag manifold $M_c$ such that $M$ is isometrically imbedded in $M_c$. If $M$ is a symmetric $R$-space, then $M_c$ is a hermitian symmetric space and the isometric imbedding is totally geodesic. If $M$ is already a complex flag manifold, then $M_c = M$.
\end{proposition}

  Using this fact that every real flag manifold $M$ can be isometrically embedded into a complex flag manifold $M_{\mathbb C}$  (in fact, it is a connected component of an antisymplectic involution on a complex flag manifold), C. U. S\'anchez defined in \cite{Sa97}   the {\it index number} $\#_{I}M$ of a real flag manifold $M$ as the maximal possible cardinality of the $p$-sets $A_{p}M$ (for a prime $p$), defined in terms of fixed points of symmetries of the complex flag manifolds restricted to the real one.   
  
As an extension of Theorem \ref{T:6.4} of M. Takeuchi, S\'anchez proved the following result in \cite{Sa97}.
 
 \begin{theorem}  Let $M$ be a real flag manifold. Then $\#_{I}M=b(M,\mathbb Z_2)$, where $b(M,\mathbb Z_2)=\sum_{i\geq 0} \dim b_i(M,\mathbb Z_2)$.
 \end{theorem}

\subsection{$k$-numbers for complex flag manifolds}\label{section 11.3}

Now, we explain very briefly the notion of $k$-number for complex flag manifold introduced  by S\'anchez in \cite{Sa93}.

 For every complex flag manifold $M_{\mathbb C}$ there exists a positive integer $k_{0}=k_{o}(M_{\mathbb C})\geq 2$ such that for each integer $k\geq k_{0}$ there exists a {\it $k$-symmetric structure} on $M_{\mathbb C}$, i.e., for each point $x\in M_{\mathbb C}$ there exists an isometry $\theta_x$ satisfying $\theta_x^k=id$ with $x$ as an isolated fixed point. 
 
 Following \cite{Ko80}, a $k$-symmetric structure is called {\it regular} it  satisfies 
 \begin{align}\label{11.1} \theta_x\circ \theta_y=\theta_z\circ \theta_z\end{align} with $z=\theta_x(y)$.
 
As an extension of two-numbers for symmetric spaces,  C. U. S\'anchez defined the {\it $k$-number} $\#_{k}(M_{\mathbb C})$ of the complex flag manifold $M_{\mathbb C}$ as the maximal possible cardinality of the so-called $k$-sets $A_{k}\subset M_{\mathbb C}$ with the property that for each $x\in A_{k}$ the corresponding $k$-symmetry fixes every point of $A_{k}$. 

S\'anchez also proved the following two results in \cite{Sa93}.

\begin{theorem} \label{T:11.2} For each complex flag manifold $M_{\mathbb C}$, one has
\begin{align}\label{11.2} \#_{k}(M_{\mathbb C})=\dim H^{*}(M_{\mathbb C},\mathbb Z_{2}).\end{align}
\end{theorem}

Let $G$ be a compact connected semisimple Lie group. Then the homogeneous spaces one obtains as orbits of $G$ under the adjoint representation on the Lie algebra of $G$ are also called {\it generalized flag manifolds}. 

It is well-known that these manifolds admit K\"ahler-Einstein metrics.
It is also well-known that every generalized flag manifold admits $k$-symmetric structure.

In a similar way, S\'anchez also proved the following.

\begin{theorem} \label{T:11.3} If $M$ is a  generalized flag manifold, then any of its $k$-symmetric structure on $M$ satisfies $\#_{k}(M)=\chi(M).$
\end{theorem}

\subsection{$k$-number and $k$-symmetric submanifolds}\label{section 11.4}

Let $\phi: M\to {\mathbb E}^m$ be an isometric embedding of a compact  Riemannian manifold into the Euclidean $m$-space.  Denote the  normal space of $M$ in $\mathbb E^m$ at $x$ by $T^\perp_x M$.
   
     If $\phi$ satisfies the following three properties: 
 \begin{itemize}
  \item[{\rm (a)}]   For each $x\in M$, there is an isometry $\sigma_x: {\mathbb E}^m\to {\mathbb E}^m$ such that 
  $\sigma_x^k =id_M$, $\sigma_x(x)=x$, and $\sigma_x |_{T^\perp_x M}$ = identity on $T^\perp_x M$; 
   
  \item[{\rm (b)}]  $\sigma_x(M)\subset M$;  and 
   
   \item[{\rm (c)}]  Let $\theta_x=\sigma_x |_M$. The collection $\{\theta_x, \,x\in M\}$ defines on $M$ a 
    Riemannian regular $s$-structure of order $k$,
      \end{itemize}
       then $M$ is called  an {\it  extrinsic $k$-symmetric submanifold} in \cite{Fe80}.
 
S\'anchez proved in \cite{Sa93} the following.
 
\begin{theorem} \label{T:11.4} If $M\subset \mathbb E^m$ is an {\it extrinsic $k$-symmetric submanifold}, then  \begin{align} \#_{k}(M)=b(M,\mathbb Z_p)\end{align}
 for any  prime number $p\geq 2$ which divides $k$.
 \end{theorem}

 \subsection{Index numbers and CW complex structures}\label{11.5} 
  A CW complex is made of basic building blocks called cells.  An $n$-dimensional closed cell is the image of an $n$-dimensional closed ball under an attaching map.  
    An $n$-dimensional open cell is a topological space that is homeomorphic to the $n$-dimensional open ball.  Closure-finite means that each closed cell is covered by a finite union of open cells. The C stands for ``closure-finite'' and the W for ``weak topology''.
 \vskip.1in

Professor Nagano and I made the following conjecture on two-numbers which was presented the first time in my report \cite[page 53]{C87}.
 
 \vskip.05in
 \noindent {\bf Conjecture. } {\it  The two-number $\#_{2}M$ is equal to the smallest number of cells that are needed  
   for a CW complex structure on a connected compact symmetric space $M$.}
  \vskip.1in

Using the convexity Theorems of Atiyah \cite{At82} and of Guillemin-Sternberg \cite{GS82} for symplectic manifolds with a hamiltonian torus action and the generalization of Duistermaat \cite{Du83} for fixed point set of antisymplectic involutions, J. Berndt, S. Console and A. Fino gave  in \cite{BCF01} alternative proofs of S\'anchez's results mentioned above.
Also, related to our conjecture given above, J. Berndt, S. Console and A. Fino proved the following result in \cite{BCF01} (see also \cite{Co}).
 
 \begin{theorem} \label{T:11.5} The index number $\#_{I}M$ is equal to the smallest number of cells 
that are needed for a CW complex structure on $M$  for each real flag manifold $M$.
 \end{theorem}

\section{Conjectures}\label{section 12}

Let $M$ be a connected compact  symmetric space.  Professor Nagano 
 and I made  the following three conjectures.

\begin{conjecture}  $\chi(M)\equiv \#_2(M)$ {\rm (mod 2)} for every connected compact  symmetric space M.
\end{conjecture}

\begin{conjecture}  $\#_{2}M=\dim H(M,\mathbb Z_{2})$  for every connected compact  symmetric space M.
\end{conjecture}

\begin{remark} M. Takeuchi  proved that Conjecture 2 is true for all symmetric $R$-spaces in \cite{Ta89}. 
On the other hand, for every compact symmetric space we could check both Conjectures 1 and 2 are correct.
\end{remark}

\begin{conjecture} $\#_{2}M =$  the smallest number of cells that are needed for a CW complex structure  for every connected compact  symmetric space M.
\end{conjecture} 

Direct computations show that  Conjecture 3 is also true for sphere, real projective space as well as for hermitian symmetric spaces.

\begin{remark} These 3 conjectures  remain open till now.\end{remark}

\section{Open problems}\label{section 13}

The total Betti numbers, $b(M;{\bf R})=\sum_{i\geq 0} b_i(M;{\bf R})$, of a simply-connected compact symmetric space $M$ satisfies (see \cite[page 53-54]{C87}):
\begin{align}\#_2M\geq b(M;{\bf R}).\end{align}

Professor T. Nagano asked   the following.
 
\begin{problem}  $b(M;{\bf R})< \#_2M$  $ \Longrightarrow$   $M$ has 2-torsion?
 \end{problem}
 
 A. Borel and J. P. Serre defined the  $p$-rank ($p$ prime), $r_{p}(G)$, of a compact  Lie group $G$ as the largest integer $h$ such that $G$ contains the direct product of $h$ cyclic groups of order $p$.   Very little were known for $p$-rank. However, Borel and Serre in \cite{BS53}  proved the next two results.
 
 \vskip.05in
 (a) $r_p(G)\geq rk(G)$ for a connected compact Lie group $G$; 

(b) $r_p(G)> rk(G)$ implies that $G$ has  $p$-torsion; 
 
 \vskip.05in
 A. M. Cohen and G. M. Seita proved, via algebraic method in \cite{CS87}, that for an algebraic group $G$ of exceptional type $G_2,F_4,E_6,E_7$, or $E_8$ over an algebraically closed field, the maximal elementary abelian $p$-subgroup $E$ of $G$ with $p\ne 2$ and of maximum rank (i.e., $rk(E)=rk(G)$) is {\it unique up to conjugacy} and of types $G_2,F_4, E_7^*$, and $E_8$ with $rk(E)=rk (G)+1$.  
 
 Besides these, very little were known for $p$-ranks of compact Lie groups. 

In view of  \cite{CS87} and \cite{Sa93}, we asked the following.
 
 \begin{problem} Is it possible to establish a geometric theory to determine 
   $p$-rank $r_{p}(G)$ for compact Lie groups for prime $p>2$?
 \end{problem}

\begin{remark} We did this for $p=2$ in \cite{CN88}.
\end{remark}

A {\it Chevalley group} is a  linear algebraic group of Lie type over finite fields related to a semisimple complex Lie algebra.   Such groups were introduced by C. Chevalley in  \cite{Ch55}. 
 For the Lie algebras $A_n, B_n, C_n, D_n$ this gave well known classical groups, but Chevalley's construction also gave groups associated to the exceptional Lie algebras $E_6, E_7, E_8, F_4$, and $G_2$. 
  The ones of type $E_6$ and of type $G_2$ were  constructed by L. E. Dickson  \cite{Di01} in 1901 and  \cite{Di05} in 1905, respectively, known as the Dickson groups. 
  
 \vskip.05in
While we were working on two-numbers during the 1980s, Professor Nagano asked the following question.
 
 \begin{problem}  Is it possible to extend our study of 2-numbers to 
 Chevalley groups? \end{problem} 

Because I don't know much about finite group theory, I don't know how to answer this question. 
I do hope that someones were able to solve this problem.

\end{document}